\documentclass[11pt]{article}

\usepackage{latexsym,color,graphicx,amsmath,amsthm,amssymb}
\usepackage{lineno}

\def\reals{{\mathbb R}}
\def\cplx{{\mathbb C}}

\def\C{{\cal C}}
\def\F{{\cal F}}

\def\XX{{1.7}}

\def\P{{\mathbb P}}
\def\O{{\mathcal O}}

\def\eps{{\varepsilon}}

\def\P{{\mathbb P}}
\newcommand{\RR}{\ensuremath{\mathbb R}}

\theoremstyle{plain}
\newtheorem{theorem}{Theorem}[section]
\newtheorem{corollary}[theorem]{Corollary}

\newtheorem{lemma}[theorem]{Lemma}

\newcommand{\Deg}{D} 
\makeatletter
\newcommand{\ProofEndBox}{{\ifhmode\unskip\nobreak\hfil\penalty50 \else
          \leavevmode\fi\quad\vadjust{}\nobreak\hfill$\Box$
            \finalhyphendemerits=0 \par}}
\makeatother

\newcommand{\proofend}{\ProofEndBox\smallskip}
\newcommand{\ignore}[1]{}

\begin{document}

\title{Incidences with curves in three dimensions\thanks{%
Work on this paper by Noam Solomon and Micha Sharir was supported by
Grants 892/13 and 260/18 from the Israel Science Foundation. Work by Micha
Sharir was also supported by Grant 2012/229 from the U.S.--Israel
Binational Science Foundation, by Grant G-1367-407.6/2016 from the 
German-Israeli Foundation for Scientific Research and Development,
by the Israeli Centers of Research
Excellence (I-CORE) program (Center No.~4/11), by the Blavatnik
Research Fund in Computer Science at Tel Aviv University and by the
Hermann Minkowski-MINERVA Center for Geometry at Tel Aviv
University. A preliminary version of an expanded version of the paper has appeared 
in {\it Proc. 28th ACM-SIAM Symposium on Discrete Algorithms} (2017), 2456--2475.}}

\author{
Micha Sharir\thanks{%
School of Computer Science, Tel Aviv University, Tel Aviv 69978,
Israel. {\sl michas@tau.ac.il} } 
\and
Noam Solomon\thanks{%
School of Computer Science, Tel Aviv University, Tel Aviv 69978,
Israel. {\sl noam.solom@gmail.com} } }

\maketitle

\begin{abstract}
We study incidence problems involving points and curves
in $\reals^3$. The current (and in fact only viable)
approach to such problems, pioneered by Guth and Katz~\cite{GK,GK2},
requires a variety of tools from algebraic geometry, most notably
(i) the polynomial partitioning technique, and (ii) the study of algebraic surfaces
that are ruled by lines or, in more recent studies~\cite{GZ}, by algebraic curves
of some constant degree. By exploiting and refining these tools, we obtain
new and improved bounds for point-curve incidence problems in $\reals^3$.

Incidences of this kind have been considered in several previous studies, 
starting with Guth and Katz's work on points and lines~\cite{GK2}. 
Our results, which are based on the work of Guth and Zahl~\cite{GZ} 
concerning surfaces that are doubly ruled by curves, provide a grand 
generalization of most of the previous results. We reconstruct the bound for 
points and lines, and improve, in certain signifcant ways, recent bounds involving
points and circles (in~\cite{SSZ}), and points and arbitrary constant-degree algebraic
curves (in~\cite{SSS}). While in these latter instances the bounds are not known (and
are strongly suspected not) to be tight, our bounds are, in a certain sense, the best
that can be obtained with this approach, given the current state of knowledge.

As an application of our point-curve incidence bound, we show that
the number of triangles spanned by a set of $n$ points in $\reals^3$ and similar 
to a given triangle is $O(n^{15/7})$, which improves the bound of Agarwal et al.~\cite{AAPS}.
Our results are also related to a study by Guth et al.~(work in progress), and have been
recently applied in Sharir et al.~\cite{ShZl} to related incidence problems in three dimensions.
\end{abstract}

\noindent {\bf Keywords.} Combinatorial geometry, incidences, the
polynomial method, infinitely ruled surfaces, algebraic geometry.

\section{Introduction} \label{sec:intro}


\subsection{The main results}

The paper studies incidences between points and constant-degree algebraic curves
in $\reals^3$. It derives several results that yield improved bounds,
and have several significant additional advantages over previous work.

The introduction states the various new results, discusses the relevant background, 
and introduces and discusses various parameters and constructs that control the 
sharpness of the derived bounds. To help the reader, we list some of the main 
results right upfront, somewhat prematurely, and supplement them with informal 
and brief explanation of the setups that we consider. Full and rigorous details
are given later in the introduction. (The numbering of the theorems below is 
their numbering in the main body of the introduction.)

\paragraph{Incidences with algebraic curves: An informal preview.}
We consider infinite families $\C_0$ of algebraic curves in $\reals^3$ of 
constant degree $E$, and assume that the curves in $\C_0$ have 
\emph{$k$ degrees of freedom}---any $k$ points determine only
$O(1)$ curves from $\C_0$ incident to all of them, and any pair 
of curves intersect at only $O(1)$ points. 

An algebraic surface $V$ is \emph{infinitely ruled} by curves from $\C_0$ 
if each \emph{generic}\footnote{%
  Genericity is a standard notion in algebraic geometry; see, e.g., 
  Cox et al.~\cite[Definition 3.6]{CLO} and also \cite[Section 2.1]{SS4d}.}
point $p\in V$ is incident to infinitely many curves 
from $\C_0$ that are fully contained in $V$. We then have:

\noindent{\bf Theorem 1.4.}
{\it
Let $P$ be a set of $m$ points and $C$ a set of $n$ irreducible
algebraic curves of constant degree $E$ in $\reals^3$, taken 
from a family $\C_0$ with $k$ degrees of freedom,
such that no surface that is infinitely ruled by curves of $\C_0$
contains more than $q$ curves of $C$, for a parameter $q<n$. Then
$$
I(P,C) = O\left(m^{\frac{k}{3k-2}}n^{\frac{3k-3}{3k-2}} +
m^{\frac{k}{2k-1}}n^{\frac{k-1}{2k-1}}q^{\frac{k-1}{2k-1}}+m+n\right),
$$
where the constant of proportionality depends on $k$, $E$, and
the complexity of the family $\C_0$.
}

See later for a precise definition of the complexity of a family $\C_0$.

We say that $\C_0$ is a family of \emph{reduced dimension} $s$ if, for each 
surface $V$ that is infinitely ruled by curves of $\C_0$, the subfamily
of the curves of $\C_0$ that are fully contained in $V$ is $s$-dimensional;
loosely, and not as general as the precise definition given later,
this means that each curve in this subfamily can be specified
algebraically in terms of $s$ real parameters. For example, the `full' dimension
of the family of circles in $\reals^3$ is $6$, but the reduced dimension 
of the subfamily of circles that lie on some fixed plane or sphere (which are 
the only surfaces that are infinitely ruled by circles~\cite{Lub}) is only $3$.
In this case we obtain the following variant of Theorem~\ref{th:mainInc3}.

\noindent{\bf Theorem 1.5.}
{\it Let $P$ be a set of $m$ points and $C$ a set of $n$ irreducible
algebraic curves of constant degree $E$ in $\reals^3$, taken from a family 
$\C_0$ that has $k$ degrees of freedom and reduced dimension $s$,
such that no surface that is infinitely ruled by curves of $\C_0$ 
contains more than $q$ of the curves of $C$. Then
$$
I(P,C) = O\left(m^{\frac{k}{3k-2}}n^{\frac{3k-3}{3k-2}}
\right) + O_\eps\left(
m^{2/3}n^{1/3}q^{1/3} +
m^{\frac{2s}{5s-4}}n^{\frac{3s-4}{5s-4}}q^{\frac{2s-2}{5s-4}+\eps}+m+n\right),
$$
for any $\eps>0$, where the first constant of proportionality depends on
$k$, $s$, $E$, and the maximum complexity of any subfamily
of $\C_0$ consisting of curves that are fully contained in some surface
that is infinitely ruled by curves of $\C_0$, and the second constant also depends on $\eps$.
}

\paragraph{Incidences with circles: An informal preview.}
Theorem 1.4 yields the fundamental incidence bound of Guth and Katz~\cite{GK2},
reviewed later in the introduction, for points and lines in $\reals^3$, 
because lines have $k=2$ degrees of freedom, and the only surfaces that 
are infinitely ruled by lines are planes; see below for more details. Theorem 1.5 
yields an improved and refined bound for incidences between points and circles 
in $\reals^3$, as compared with the previous result of \cite{SSZ} (see also
another bound in \cite{Zahl}, which will be discussed later in the introduction):

\noindent
{\bf Theorem 1.6.} 
{\it
Let $P$ be a set of $m$ points and $C$ a set of $n$ circles in $\reals^3$, so that
no plane or sphere contains more than $q$ circles of $C$. Then
$$
I(P,C) = O\left( m^{3/7}n^{6/7} + m^{2/3}n^{1/3}q^{1/3} +
m^{6/11}n^{5/11}q^{4/11}\log^{2/11}(m^3/q) + m + n \right) .
$$
}

We now proceed to the full introduction, in which these results, and several others, will be given,
along with the background and detailed discussion of the various parameters and notions mentioned above.

A preliminary version of the paper has appeared as part of the conference paper~\cite{SS:soda}.


\paragraph{The setup: Incidences between points and curves in three dimensions.}
Let $P$ be a set of $m$ points and $C$ a set of $n$ irreducible algebraic
curves of constant degree in $\reals^3$. We consider the problem of obtaining
sharp incidence bounds between the points of $P$ and the curves of $C$.
This is a major topic in incidence geometry since the groundbreaking work of
Guth and Katz~\cite{GK2} on point-line incidences in $\reals^3$,
with many follow-up studies, some of which are reviewed below.
Building on the recent work of Guth and Zahl~\cite{GZ}, which bounds
the number of \emph{2-rich points} determined by a set of bounded-degree algebraic
curves in $\reals^3$ (i.e., points incident to at least two of the given curves),
we are able to generalize Guth and Katz's point-line incidence bound to a general
bound on the number of incidences between points and bounded-degree irreducible
algebraic curves that satisfy certain natural assumptions, discussed in detail below.

\subsection{Background}

\paragraph{Points and curves, the planar case.}
The case of incidences between points and curves has a rich history,
starting with the case of points and lines in the
plane, studied in the seminal paper of Szemer\'edi and Trotter~\cite{ST}, 
and later also in~\cite{CEGSW,Sze}, where the worst-case tight bound on the
number of incidences is $\Theta(m^{2/3}n^{2/3}+m+n)$, where $m$ is
the number of points and $n$ is the number of lines. Still in the
plane, Pach and Sharir~\cite{PS} extended this bound to incidence
bounds between points and curves with $k$ \emph{degrees of freedom}.
These are curves with the property that, for each set of $k$ points,
there are only $\mu=O(1)$ curves that pass through all of them, and
each pair of curves intersect in at most $\mu$ points; $\mu$ is
called the \emph{multiplicity} (of the degrees of freedom).
\begin{theorem}[Pach and Sharir \cite{PS}] \label{th:PS}
Let $P$ be a set of $m$ points in $\reals^2$ and let $C$ be a set
of $n$ bounded-degree algebraic curves in $\reals^2$ with $k$
degrees of freedom and with multiplicity $\mu$. Then
$$
I(P,C) = O\left(m^{\frac{k}{2k-1}}n^{\frac{2k-2}{2k-1}}+m+n\right)
,
$$
where the constant of proportionality depends on $k$ and $\mu$.
\end{theorem}

\noindent{\bf Remark.} The result of Pach and Sharir holds for more
general families of curves, not necessarily algebraic, but, since
algebraicity will be assumed in three dimensions, we assume it also
in the plane.

Note that the Szemer\'edi--Trotter bound is a special case of 
Theorem~\ref{th:PS}, since lines have two degrees of freedom. However,
except for the case of lines (and $k=2$), the bound is not known, and 
is strongly suspected not, to be tight in the worst case. Indeed, 
in a series of papers during the 2000's~\cite{ANPPSS, ArS, MT}, 
an improved bound has been obtained
for incidences with circles, parabolas, or other families of curves
with certain properties (see~\cite{ANPPSS} for the precise
formulation). Specifically, for a set $P$ of $m$ points and a set
$C$ of $n$ circles, or parabolas, or similar curves, we have
\begin {equation}
\label {eq:impbo} I(P,C) =
O(m^{2/3}n^{2/3}+m^{6/11}n^{9/11}\log^{2/11}(m^3/n)+m+n).
\end {equation}
Some further (slightly) improved bounds, over the bound in
Theorem~\ref{th:PS}, for more general families of curves in the
plane have been obtained by Chan~\cite{Chan, Chan2} and by
Bien~\cite{Bien}. They are, however, considerably weaker than the
bound in (\ref{eq:impbo}).

Recently, Sharir and Zahl~\cite{SZ} have considered general families
of constant-degree algebraic curves in the plane that belong to an
\emph{$s$-dimensional family of curves}. Extending the definition
given above, this means that each curve in that family
can be represented by a constant number of real parameters, so that,
in this parametric space, the points representing the curves lie in
an $s$-dimensional algebraic variety $\F$ of some constant degree
(to which we refer as the ``complexity'' of $\F$).
See~\cite{SZ} for more details.
\begin{theorem}[Sharir and Zahl~\cite{SZ}]\label{incPtCu}
Let $C$ be a set of $n$ algebraic plane curves that belong to
an $s$-dimensional family $\F$ of curves of maximum constant degree
$E$, no two of which share a common irreducible component, and let
$P$ be a set of $m$ points in the plane. Then, for any $\eps>0$,
the number $I(P,C)$ of incidences between the points of
$P$ and the curves of $C$ satisfies
\begin{equation*} \label{newIncBd}
 I(P,C) = O\Big(m^{\frac{2s}{5s-4}} n^{\frac{5s-6}{5s-4}+\eps} + m^{2/3}n^{2/3} + m + n\Big) ,
\end{equation*}
where the constant of proportionality depends on $\eps$, $s$, $E$,
and the complexity of the family $\F$.
\end{theorem}
Except for the factor $O(n^\eps)$, this is a significant improvement
over the bound in Theorem~\ref{th:PS} (for $s\ge 3$), in cases where
the assumptions in Theorem~\ref{incPtCu} imply (as they often do)
that $C$ has $k=s$ degrees of freedom. Concretely, when $k=s$, we
obtain an improvement, except for the factor $n^\eps$, which is an
artifact of the proof of Theorem~\ref{incPtCu}, for the entire
``meaningful'' range $n^{1/s} \le m \le n^2$, in which the bound is
superlinear. The factor $n^\eps$ makes the bound in \cite{SZ}
slightly weaker only when $m$ is close to the lower end $n^{1/s}$ of
that range. Note also that for circles (where $s=3$), the bound in
Theorem~\ref{incPtCu} nearly coincides with the slightly more
refined bound (\ref{eq:impbo}).

\paragraph{Incidences with curves in three dimensions.}
The seminal work of Guth and Katz~\cite{GK2} establishes the sharper
bound $O(m^{1/2}n^{3/4} + m^{2/3}n^{1/3}q^{1/3} + m + n)$ on the
number of incidences between $m$ points and $n$ lines in $\reals^3$,
provided that no plane contains more than $q$ of the given lines.\footnote{%
  This bound is not explicitly mentioned in \cite{GK2}, but it
  immediately follows from the analysis given there.}
This has lead to many recent works on incidences between points and
lines or other curves in three and higher dimensions;
see~\cite{CPS,GZ,SSS,SSZ,SS4d,SS4dv,Zahl} for a sample of these results.
Most relevant to our present study are the works of Sharir, Sheffer,
and Solomon~\cite{SSS} on incidences between points and curves in
any dimension, the work of Sharir, Sheffer, and Zahl~\cite{SSZ} on
incidences between points and circles in three dimensions, the
work of Sharir and Solomon~\cite{SS4d} on incidences between points
and lines in four dimensions, and the recent work of Do and Sheffer~\cite{DS}
that presents more general bounds for incidences with curves and surfaces
in any dimension, as well as several other studies of
point-line incidences by the authors~\cite{SS3d, SS4dv}.

Of particular significance is the recent work of Guth and Zahl~\cite{GZ} 
on the number of 2-rich points in a collection of
curves, namely, points incident to at least two of the given curves.
For the case of lines, Guth and Katz~\cite{GK2} have shown that the
number of such points is $O(n^{3/2})$, when no plane or regulus
contains more than $O(n^{1/2})$ lines. Guth and Zahl obtain the same
asymptotic bound for general algebraic curves, under analogous (but
stricter) restrictive assumptions.

The new bounds that we will derive require the extension to three dimensions of the notions
of having $k$ degrees of freedom and of being an $s$-dimensional family of curves.
The definitions of these concepts, as given above for the planar case, extend,
more or less verbatim, to three (or higher) dimensions, but, even in typical situations, 
these two concepts do not coincide anymore. For example, lines in three dimensions have 
two degrees of freedom, but they form a $4$-dimensional family of curves (this is the 
number of parameters needed to specify a line in $\reals^3$).

\subsection{Our results} \label{sec:res}

Before we state our results, we discuss three
notions that are used in these statements. These are the notions of
\emph{$k$ degrees of freedom} (already mentioned above), of
\emph{constructibility}, and of \emph{surfaces infinitely ruled by curves.}

\paragraph{$k$ degrees of freedom.}
Let $\C_0$ be an infinite family of irreducible algebraic curves of constant
degree $E$ in $\RR^3$. Formally, in complete analogy with the planar case, we say
that $\C_0$ has $k$ \emph{degrees of freedom} with \emph{multiplicity} $\mu$,
where $k$ and $\mu$ are constants, if (i) for every tuple of $k$ points in $\RR^3$
there are at most $\mu$ curves of $\C_0$ that are incident to all $k$ points, and
(ii) every pair of curves of $\C_0$ intersect in at most $\mu$ points. As in~\cite{PS},
the bounds that we derive depend more significantly on $k$ than on $\mu$---see below.

\paragraph{Constructibility.}
In the statements of the following theorems, we also assume that $\C_0$
is a \emph{constructible} family of curves. This notion generalizes the notion 
of being algebraic, and is discussed in detail in Guth and Zahl~\cite{GZ}. 
Informally, a set $Y\subset \cplx^d$ is constructible
if it is a Boolean combination of algebraic sets. The formal definition goes as
follows (see, e.g., Harris~\cite[Lecture 3]{Harris}).
For $z\in\cplx$, define $v(0)=0$ and $v(z)=1$ for $z\ne 0$. Then $Y\subseteq \cplx^d$,
for some fixed $d$, is a \emph{constructible set} if there exist an integer $J_Y$,
a finite set of polynomials
$f_j : \cplx^d \to \cplx$, for $j = 1,\ldots,J_Y$, and a subset $B_Y \subset \{0,1\}^{J_Y}$,
so that $x \in Y$ if and only if $\left(v(f_1(x)),\ldots, v(f_{J_Y}(x))\right) \in B_Y$.

When we apply this definition to a set of curves, we think of them as points in some
parametric (complex) $d$-space, where $d$ is the number of parameters needed to specify
a curve. When $J_Y=1$ we get all the algebraic hypersurfaces (that admit the implied
$d$-dimensional representation) and their complements. An $s$-dimensional family of curves,
for $s<d$, is obtained by taking $J_Y=d-s$ and $B_Y=\{0\}^{J_Y}$. In doing so, we assume that
the curves that we obtain are complete intersections. Following Guth and Zahl (see also a 
comment to that effect in the appendix), this involves no loss of generality, because every 
curve is contained in a curve that is a complete intersection. 
In what follows, when we talk about constructible sets, we implicitly assume 
that the ambient dimension $d$ is constant.

The constructible sets form a Boolean algebra. This means that finite unions and
intersections of constructible sets are constructible, and the complement of a
constructible set is constructible.
Another fundamental property of constructible sets is that, over $\cplx$, the projection
of a constructible set is constructible; this is known as Chevalley's theorem
(see Harris~\cite[Theorem 3.16]{Harris} and Guth and Zahl~\cite[Theorem 2.3]{GZ}).
If $Y$ is a constructible set, we define the \emph{complexity} of $Y$ to be 
$\min(\deg f_1 + \cdots + \deg f_{J_Y})$, where the minimum
is taken over all representations of $Y$, as described above. As just observed,
constructibility of a family $\C_0$ of curves extends the notion of $\C_0$ being
$s$-dimensional. One of the main motivations for using the notion of constructible
sets (rather than just $s$-dimensionality) is the fact, established by
Guth and Zahl~\cite[Proposition 3.3]{GZ}, that the set $\C_{3,E}$ of \emph{irreducible}
curves of degree at most $E$ in complex 3-dimensional space (either affine or projective)
is a constructible set of constant complexity that depends only on $E$. Moreover,
Theorem~\ref{th:ruledgz}, one of the central technical tools that we use in our analysis
(see below for its statement and the appendix for its proof), 
holds for constructible families of curves.

\smallskip

\paragraph{The connection between degrees of freedom and constructibility/dimensionality.}
Loosely speaking, in the plane the number of degrees of freedom and the dimensionality 
of a family of curves tend to be equal. In three dimensions the situation is different. 
This is because the constraint that a curve $\gamma$ passes through a point $p$ imposes 
two equations on the parameters defining $\gamma$. We therefore expect the number of 
degrees of freedom to be half the dimensionality. A few instances that illustrate
this connection are:

\medskip
\noindent
(i) Lines in three dimensions have two degrees of freedom, and they form a
$4$-dimensional family of curves (this is the number of parameters needed to
specify a line in $\reals^3$).

\medskip
\noindent
(ii) Circles in three dimensions have three degrees of
freedom, and they form a $6$-dimensional family of curves (e.g., one needs three
parameters to specify the plane containing the circle, two additional parameters to
specify its center, and a sixth parameter for its radius). 

\medskip
\noindent
(iii) Conic sections have five degrees of freedom, but they form an $8$-dimensional 
family of curves, as is easily checked (we need three parameters for the plane
containing the curve, and five parameters to define a conic section in that plane). 
This discrepancy (for conic sections) is explained by noting that four points are 
not sufficient to define the curve, because the first three determine the plane 
containing it, so the fourth point, if at all coplanar with the first three, only 
imposes one constraint on the parameters of the curve.

\medskip
\noindent{\bf Remark.} The definition of constructibility is given
over the complex field $\cplx$. This is in accordance with most of
the basic algebraic geometry tools, which have been developed over
the complex field. Some care has to be exercised when applying them
over the reals. For example, Theorem~\ref{th:ruledgz}, one of the
central technical tools that we use in our analysis, as well as the
related results of Guth and Zahl~\cite{GZ}, apply over the complex field,
but not over the reals. On the other hand, when we apply the
partitioning method of~\cite{GK2} (as in the proof of
Theorems~\ref{th:mainInc3}) or when we use
Theorem~\ref{incPtCu}, we (have to) work over the reals.

It is a fairly standard practice in algebraic geometry that handles
a real algebraic variety $V$, defined by real polynomials, by
considering its complex counterpart $V_\cplx$, namely the set of
complex points at which the polynomials defining $V$ vanish. The
rich toolbox that complex algebraic geometry has developed allows
one to derive various properties of $V_\cplx$, which, with some care, can
usually be transported back to the real variety $V$.

This issue arises time and again in this paper. Roughly speaking, we approach
it as follows. We apply the polynomial partitioning technique to the given
sets of points and of curves or surfaces, in the original real (affine) space,
as we should. Within the cells of the partitioning we then apply some
field-independent argument, based either on induction or on some ad-hoc
combinatorial argument. Then we need to treat points that lie on the zero set
of the partitioning polynomial. We can then switch to the complex field, when
it suits our purpose, noting that this step preserves all the real incidences;
at worst, it might add additional incidences involving the non-real portions
of the variety and of the curves. Hence, the bounds that we obtain
for this restrictive setup transport, more or less verbatim, to the real case too.

\paragraph{Surfaces infinitely ruled by curves.}
Back in three dimensions, a surface $V$ is (singly, doubly, or
infinitely) ruled by some family $\Gamma$ of curves of degree at
most $E$, if every \emph{generic} point $p\in V$ is incident to (at least one, at
least two, or infinitely many) curves of $\Gamma$ that are fully
contained in $V$. The connection between ruled surface theory and
incidence geometry goes back to the pioneering work of Guth and
Katz~\cite{GK2} and shows up in many subsequent works. See Guth's
survey~\cite{Gut:surv} and book~\cite{Gut:book}, and
Koll\'ar~\cite{Kollar} for details.

In most of the previous works, only singly-ruled and doubly-ruled
surfaces have been considered. Looking at infinitely-ruled surfaces
adds a powerful ingredient to the toolbox, as will be demonstrated
in this paper.

We recall that the only surfaces that are infinitely ruled by lines
are planes (see, e.g., Fuchs and Tabachnikov~\cite[Corollary 16.2]{FT}), 
and that the only surfaces that are infinitely ruled by circles are 
spheres and planes (see, e.g., Lubbes~\cite[Theorem 3]{Lub} and 
Schicho~\cite{Sch}; see also Skopenkov and
Krasauskas~\cite{SK} for recent work on \emph{celestials}, namely
surfaces \emph{doubly} ruled by circles, and Nilov and
Skopenkov~\cite{NS13}, proving that a surface that is ruled by a
line and a circle through each (generic) point is a quadric). 
It should be noted that,
in general, for this definition to make sense, it is important to
require that the degree $E$ of the ruling curves be smaller than $\deg(V)$. 
Otherwise, every variety $V$ is infinitely ruled by,
say, the curves $V\cap h$, for hyperplanes $h$, having the same
degree as $V$. A challenging open problem is to characterize all the
surfaces that are infinitely ruled by algebraic curves of degree at
most $E$ (or by certain special classes thereof). However, the
following result of Guth and Zahl provides a useful necessary
condition for this property to hold.
\begin{theorem} [Guth and Zahl~\cite{GZ}] \label{th:gz}
Let $V$ be an irreducible surface, and suppose that it is
\emph{doubly ruled} by curves of degree at most $E$. Then $\deg(V) \le 100E^2$.
\end{theorem}
In particular, an irreducible surface that is infinitely ruled by
curves of degree at most $E$ is doubly ruled by these curves, so its
degree is at most $100E^2$. Therefore, if $V$ is irreducible of
degree $D$ larger than this bound, $V$ cannot be infinitely ruled by
curves of degree at most $E$. This leaves a gray zone, in which the
degree of $V$ is between $E$ and $100E^2$. We would like to conjecture 
that if $V$ is an irreducible variety that is infinitely ruled by 
degree-$E$ curves then its degree is $O(E)$. Being unable to establish 
this conjecture, we leave it as a challenging open problem for further research.

\paragraph{Our results.}

We can now state our main results on point-curve incidences.
\begin{theorem} \label{th:mainInc3}
Let $P$ be a set of $m$ points and $C$ a set of $n$ irreducible
algebraic curves of constant degree $E$, taken from a constructible family
$\C_0$, of constant complexity, with $k$ degrees of freedom (and some multiplicity
$\mu$) in $\RR^3$, such that no surface that is infinitely ruled by curves of $\C_0$
contains more than $q$ curves of $C$, for a parameter $q<n$. Then
\begin{equation} \label{eq:mainInc3}
I(P,C) = O\left(m^{\frac{k}{3k-2}}n^{\frac{3k-3}{3k-2}} +
m^{\frac{k}{2k-1}}n^{\frac{k-1}{2k-1}}q^{\frac{k-1}{2k-1}}+m+n\right),
\end{equation}
where the constant of proportionality depends on $k$, $\mu$, $E$, and
the complexity of the family $\C_0$.
\end{theorem}

\noindent{\bf Remarks. (1)} In certain favorable situations, such as
in the cases of lines or circles, discussed above, the surfaces that
are infinitely ruled by curves of $\C_0$ have a simple
characterization. In such cases the theorem has a stronger flavor,
as its assumption on the maximum number of curves on a surface has
to be made only for this concrete kind of surfaces. For example, as
already noted, for lines (resp., circles) we only need to require
that no plane (resp., no plane or sphere) contains more than $q$ of
the curves. In general, as mentioned, characterizing
infinitely-ruled surfaces by a specific family of curves is a
difficult task. Nevertheless, we can overcome this issue by
replacing the assumption in the theorem by a more restrictive one,
requiring that no surface that is infinitely ruled by curves of
degree at most $E$ contain more than $q$ curves of $C$. By
Theorem~\ref{th:gz}, any infinitely ruled surface of this kind must
be of degree at most $100E^2$. Hence, an even simpler (albeit
weaker) formulation of the theorem is to require that no surface of
degree at most $100E^2$ contains more than $q$ curves of $C$. This
can indeed be much weaker: In the case of circles, say, instead of
making this requirement only for planes and spheres, we now have to
make it for every surface of degree at most $400$.
See also \cite{GZ} for a similar approach.

\smallskip

\noindent{\bf (2)} In several recent works (see~\cite{Gut,SSS,SSZ}),
the assumption in the theorem is replaced by a much more restrictive
assumption, that no surface of degree at most $c_\eps$ contains more
than $q$ given curves, where $c_\eps$ is a constant that depends on
another prespecified parameter $\eps>0$ (where $\eps$ appears in the
exponents in the resulting incidence bound), and is typically very
large (and increases as $\eps$ becomes smaller). Getting rid of such
an $\eps$-dependent constant (and of the $\eps$ in the exponent) is a significant feature of
Theorem~\ref{th:mainInc3}.

\smallskip

\noindent{\bf (3)}
As already mentioned, Theorem~\ref{th:mainInc3} generalizes the incidence 
bound of Guth and Katz~\cite{GK2}, obtained for the case of lines. 
In this case, lines have $k=2$ degrees of freedom, they certainly 
form a constructible (in fact, a $4$-dimensional) family of curves,
and, as just noted, planes are the only surfaces in $\RR^3$ that are infinitely
ruled by lines. Thus, in this special case, both the assumptions and the bound in
Theorem~\ref{th:mainInc3} are identical to those in Guth and Katz~\cite{GK2}.
That is, if no plane contains more than $q$ input lines, the number of incidences is
$O(m^{1/2}n^{3/4}+m^{2/3}n^{1/3}q^{1/3}+m+n)$.

\paragraph{Improving the bound.}
The bound in Theorem~\ref{th:mainInc3} can be further improved, if
we also throw into the analysis the dimensionality $s$ of the family
$\C_0$. Actually, as will follow from the proof, the dimensionality
that will be used is only that of any subset of $\C_0$ whose members
are fully contained in some variety that is infinitely ruled by
curves of $\C_0$. As just noted, such a variety must be of constant
degree (at most $100E^2$, or smaller as in the cases of lines and
circles), and the additional constraint that the curves be contained
in the variety can typically be expected to reduce the
dimensionality of the family.

For example, if $\C_0$ is the collection of all circles in $\reals^3$, then, since
the only surfaces that are infinitely ruled by circles are spheres and planes, the
subfamily of all circles that are contained in some sphere or plane is only
$3$-dimensional (as opposed to the entire $\C_0$, which is $6$-dimensional).

We capture this setup by saying that $\C_0$ is a family of \emph{reduced dimension}
$s$ if, for each surface $V$ that is infinitely ruled by curves of $\C_0$, the subfamily
of the curves of $\C_0$ that are fully contained in $V$ is $s$-dimensional. We then
obtain the following variant of Theorem~\ref{th:mainInc3}.
\begin{theorem} \label{th:imprInc3}
Let $P$ be a set of $m$ points and $C$ a set of $n$ irreducible
algebraic curves of constant degree $E$, taken from a constructible family
$\C_0$ with $k$ degrees of freedom (and some multiplicity $\mu$) in
$\RR^3$, such that no surface that is infinitely ruled by curves of
$\C_0$ contains more than $q$ of the curves of $C$, and assume further that
$\C_0$ is of reduced dimension $s$. Then
\begin{equation} \label{eq:imprInc3}
I(P,C) = O\left(m^{\frac{k}{3k-2}}n^{\frac{3k-3}{3k-2}}
\right) + O_\eps\left(
m^{2/3}n^{1/3}q^{1/3} +
m^{\frac{2s}{5s-4}}n^{\frac{3s-4}{5s-4}}q^{\frac{2s-2}{5s-4}+\eps}+m+n\right),
\end{equation}
for any $\eps>0$, where the first constant of proportionality depends on
$k$, $\mu$, $s$, $E$, and the maximum complexity of any subfamily
of $\C_0$ consisting of curves that are fully contained in some surface
that is infinitely ruled by curves of $\C_0$, and the second constant also depends on $\eps$.
\end {theorem}

\noindent{\bf Remarks. (1)}
Theorem~\ref{th:imprInc3} is an improvement of Theorem~\ref{th:mainInc3} when
$s\le k$ and $m>n^{1/k}$, in cases where $q$ is sufficiently large so as to make
the second term in (\ref{eq:mainInc3}) dominate the first term; for smaller values
of $m$ the bound is always linear. This is true except for the term $q^\eps$, which
is an artifact of the proof of the theorem, and which
affects the bound only when $m$ is very close to $n^{1/k}$ (when $s=k$). When $s>k$
we get a threshold exponent $\beta = \frac{5s-4k-2}{ks-4k+2s}$ (which becomes $1/k$
when $s=k$), so that the bound in Theorem~\ref{th:imprInc3} is stronger (resp., weaker)
than the bound in Theorem~\ref{th:mainInc3} when $m > n^\beta$ (resp., $m < n^\beta$),
again, up to the extra factor $q^\eps$.

\smallskip

\noindent{\bf (2)}
The bounds in Theorems~\ref{th:mainInc3} and~\ref{th:imprInc3} improve, in three
dimensions, the recent result of Sharir, Sheffer, and Solomon~\cite{SSS}, in three
significant ways:
{\bf(i)} The leading terms in both bounds are essentially the same, but
our bound is sharper, in that it does not include the factor $O(n^\eps)$ appearing
in \cite{SSS}. {\bf (ii)} The assumption here, concerning the number of curves on a
low-degree surface, is much weaker than the one made in~\cite{SSS}, where
it was required that no surface of some (constant but potentially very large)
degree $c_\eps$, that depends on $\eps$, contains more than $q$ curves of $C$.
(See also Remark (2) following Theorem~\ref{th:mainInc3}.)
{\bf (iii)} The two variants of the non-leading terms here are significantly smaller
than those in~\cite{SSS}, and, in a certain sense (that will be elaborated
following the proof of Theorem~\ref{th:imprInc3}) are best possible.

\smallskip

\noindent{\bf (3)}
The results and techniques in this paper have recently been applied
in Sharir et al.~\cite{ShZl} to incidence problems in three dimensions
with curves that have `almost two degrees of freedom' (a notion defined in that paper).

\smallskip
\paragraph{Point-circle incidences in $\reals^3$.}

Theorem~\ref{th:imprInc3} yields a new bound for the case of incidences between points and circles
in $\reals^3$, which improves over the previous bound of Sharir,
Sheffer, and Zahl~\cite{SSZ}. Specifically, we have:
\begin{theorem} \label{ptcirc}
Let $P$ be a set of $m$ points and $C$ a set of $n$ circles in $\reals^3$, so that
no plane or sphere contains more than $q$ circles of $C$. Then
$$
I(P,C) = O\left( m^{3/7}n^{6/7} + m^{2/3}n^{1/3}q^{1/3} +
m^{6/11}n^{5/11}q^{4/11}\log^{2/11}(m^3/q) + m + n \right) .
$$
\end{theorem}
Note that the bound in Theorem~\ref{ptcirc} is slightly better than the
bound obtained from Theorem~\ref{th:imprInc3} for $s=3$, in that it
replaces the factor $q^\eps$ by the factor $\log^{2/11}(m^3/q)$.
This is a consequence of the proof technique and of the improved 
point-circle incidence bound in the plane (\ref{eq:impbo}).
See Section~\ref{se:sim} for details.

Here too we have the three improvements noted in Remark (2) above. 
In particular, in the sense of part (iii) of that remark, the new bound 
is, in a rather vague sense, ``best possible'' with respect to the best 
known bound (\ref{eq:impbo}) for the planar or spherical cases. 
See Section~\ref{se:sim} for details. We remark, though, that a recent
work of Zahl~\cite{Zahl} gives a different point-circle incidence bound,
namely,
\[
O^*\left( m^{1/2}n^{3/4} + m^{2/3}n^{13/15} + m^{1/3}n^{8/9}
+ nq^{2/3} + m \right) 
\]
where the notation $O^*$ hides subpolynomial factors,
which is better for a certain range of $m$, $n$, and $q$.

Theorem~\ref{ptcirc} has an interesting application to the problem of 
bounding the number of similar triangles spanned by a set of $n$ points 
in $\reals^3$. It yields the bound $O(n^{15/7})$, which improves the bound 
of Agarwal et al.~\cite{AAPS}. Again, see Section~\ref{se:sim} for details.

\subsection{The main techniques}

There are three main ingredients used in
our approach. The first ingredient, already mentioned in the context
of planar point-curve incidences, is the techniques of Pach and
Sharir~\cite{PS} (given in Theorem~\ref{th:PS}), and of Sharir and
Zahl~\cite{SZ} (Theorem~\ref{incPtCu}) concerning incidences between
points and algebraic curves in the plane. 

The second ingredient, relevant to the proof of 
Theorems~\ref{th:mainInc3} and \ref{th:imprInc3},
is the \emph{polynomial partitioning technique} of Guth and
Katz~\cite{GK2}, and its more recent extension by Guth~\cite{Guth},
which yields a divide-and-conquer mechanism via space decomposition
by the zero set of a suitable polynomial. This will produce
subproblems that will be handled recursively, and will leave us with
the overhead of analyzing the incidence pattern involving the points
that lie on the zero set itself, which in turn can be handled using
Theorem~\ref{th:ruledgz} stated below. We assume familiarity of the 
reader with these results; more details will be given in the applications 
of this technique in the proofs of the aforementioned theorems.

The third ingredient arises in the proof of Theorems~\ref{th:mainInc3} and
\ref{th:imprInc3}, where we argue that a ``generic'' point on a
variety $V$, that is not infinitely ruled by constant-degree curves of some given family,
as in the statement of the theorems, is incident to at most a
constant number of the given curves that are fully contained in $V$.
Moreover, we can also control the number and structural properties of ``non-generic'' points.

Before formally stating, in detail, the technical properties that we need, we
review a few notations.

Fix a constructible subfamily $\C_0$ of the family $\C_{3,E}$ of the 
irreducible curves of degree at most $E$ in 3-dimensional space, and 
a trivariate polynomial $f$. Following Guth and Zahl~\cite[Section 9]{GZ}, 
we call a point $p\in Z(f)$ a \emph{$(t, \C_0,
r)$-flecnode}, if there are at least $t$ curves $\gamma_1,\ldots,
\gamma_t \in \C_0$, such that, for each $i=1,\ldots, t$, (i)
$\gamma_i$ is incident to $p$, (ii) $p$ is a non-singular point of
$\gamma_i$, and (iii) $\gamma_i$ osculates to $Z(f)$ to order $r$ at
$p$. This is a generalization of the notion of a \emph{flecnodal
point}, due to Salmon~\cite[Chapter XVII, Section III]{salmon} (see
also \cite{GK2, SS4d} for more details). Our analysis requires the
following theorem. It is a consequence of the analysis of Guth and
Zahl~\cite[Corollary 10.2]{GZ}, which itself is a generalization of
the Cayley--Salmon theorem on surfaces ruled by lines (see, e.g.,
Guth and Katz~\cite{GK2}), and is closely related to
Theorem~\ref{th:gz} (also due to Guth and Zahl~\cite{GZ}). The
novelty in this theorem is that it addresses surfaces that are
\emph{infinitely ruled} by certain families of curves, whereas the
analysis in~\cite{GZ} only handles surfaces that are \emph{doubly
ruled} by such curves.
\begin{theorem} \label{th:ruledgz}
(a) For given integer parameters $c$ and $E$, there are constants
$c_1=c_1(c,E)$, $r=r(c,E)$, and $t=t(c,E)$, such that the following
holds. Let $f$ be a complex irreducible polynomial of degree\footnote{%
  See the appendix for the precise relationship between $D$ and $E$.}
$D \gg E$, and let $\C_0 \subset \C_{3,E}$ be a constructible set of
complexity at most $c$. If there exist at least $c_1 D^2$ curves of
$\C_0$, such that each of them is contained in $Z(f)$ and contains
at least $c_1 D$ points on $Z(f)$ that are $(t,\C_0,r)$-flecnodes,
then $Z(f)$ is infinitely ruled by curves from $\C_0$.

\noindent{(b)} In particular, if $Z(f)$ is not infinitely ruled by curves from
$\C_0$ then, except for at most $c_1 D^2$ \emph{exceptional} curves,
every curve in $\C_0$ that is fully contained in $Z(f)$ is incident
to at most $c_1 D$ $t$-rich points, namely points that are incident 
to at least $t$ curves in $\C_0$ that are also fully contained in $Z(f)$.
\end{theorem}
Note that (b) follows from (a) because, by definition, a $t$-rich point is
a $(t,\C_0,r)$-flecnode, for any $r\ge 1$ (a curve fully contained in the variety
osculates to it to any order).

Note that, by making $c_1$ sufficiently large (specifically, choosing $c_1 > E$),
the assumption that each of the $c_1D^2$ curves in the premises of the theorem is
fully contained in $Z(f)$ follows (by B\'ezout's
theorem) from the fact that each of them contains at least
$c_1D$ points on $Z(f)$.
Although the theorem is a (heretofore unstated) corollary of the work of Guth and Zahl in~\cite{GZ}, 
we review (in the appendix) the machinery needed for its proof, and sketch a brief version
of the proof itself, for the convenience of the reader and in the interest of completeness.

\section{Proofs of Theorems~\protect{\ref{th:mainInc3}} and~\protect{\ref{th:imprInc3}}}
\label{sec:curves}

The proofs of both theorems are almost identical, and they differ in
only one step in the analysis. We will give a full proof of
Theorem~\ref{th:mainInc3}, and then comment on the few modifications
that are needed to establish Theorem~\ref{th:imprInc3}.

\paragraph{Proof of Theorem~\ref{th:mainInc3}.}
Since the family $C$ has $k$ degrees of freedom with multiplicity
$\mu$, the incidence graph $G(P,C)$, as a subgraph of $P\times C$,
does not contain $K_{k,\mu+1}$ as a subgraph. The K{\H
o}v\'ari-S\'os-Tur\'an theorem (e.g., see~\cite[Section 4.5]{Mat02})
then implies that $I(P,C) = O(mn^{1-1/k}+n)$, where the constant of
proportionality depends on $k$ (and $\mu$). We refer to this as the
\emph{naive bound} on $I(P,C)$. In particular, when $m=O(n^{1/k})$,
we get $I(P,C) = O(n)$. We may thus assume that $m\ge a'n^{1/k}$,
for some absolute constant $a'$.

The proof proceeds by double induction on $n$ and $m$, and establishes the bound
\begin{equation} \label{eqind}
I(P,C) \le A \left(m^{\frac{k}{3k-2}}n^{\frac{3k-3}{3k-2}} +
m^{\frac{k}{2k-1}}n^{\frac{k-1}{2k-1}}q^{\frac{k-1}{2k-1}}+m+n\right) ,
\end{equation}
for a suitable constant $A$ that depends on $k$, $\mu$, $E$, and the complexity of $\C_0$.

The base case for the outer induction on $n$ is $n\le n_0$, for a
suitable sufficiently large constant threshold $n_0$ that will be set later.
The bound (\ref{eqind}) clearly holds in this case if we choose $A \ge n_0$.

The base case for the inner induction on $m$ is $m\le a'n^{1/k}$,
in which case the naive bound implies that $I(P,C) = O(n)$,
so (\ref{eqind}) holds with a sufficiently large choice of $A$.
Assume then that the bound (\ref{eqind}) holds for all sets $P'$,
$C'$ with $|C'|<n$ or with $|C'|=n$ and $|P'|<m$, and let $P$ and
$C$ be sets of sizes $|P|=m$, $|C|=n$, such that $n>n_0$, and $m>a'n^{1/k}$.

It is instructive to notice that the two terms
$m^{\frac{k}{3k-2}}n^{\frac{3k-3}{3k-2}}$ and $m$ in (\ref{eqind})
compete for dominance; the former (resp., latter) dominates when
$m\le n^{3/2}$ (resp., $m\ge n^{3/2}$). One therefore has to treat
these two cases somewhat differently; see below and also in earlier
works~\cite{GK2,SS3d}.

\smallskip
\noindent {\bf Applying the polynomial partitioning technique.}
We construct a \emph{partitioning polynomial} $f$ for the set $C$ of curves,
as in the recent variant of the polynomial partitioning technique, due
to Guth~\cite{Guth}. Specifically, we choose a degree

\begin{equation}
\label{eq:pocu} D = \begin{cases} c m^{\frac {k} {3k-2}}/n^{\frac 1
{3k-2}} , & \text{for $a'n^{1/k}\le m\le an^{3/2}$} ,
\\ cn^{1/2} , & \text{for $m > an^{3/2}$} ,
\end{cases}
\end {equation}
for suitable constants $c$, $a$, and the previously introduced constant
$a'$ (all of whose values will be set later), and obtain a polynomial $f$ of 
degree at most $D$, such that each of the $O(D^3)$ (open) connected components 
of $\reals^3\setminus Z(f)$ is crossed by at most $O(n/D^2)$ curves of
$C$, where the former constant of proportionality is absolute, and
the latter one depends on $E$. Note that in both cases $1\le D\ll
n^{1/2}$, if $a$, $a'$, and $c$ are chosen appropriately.  Denote
the cells of the partition as $\tau_1, \ldots, \tau_u$, for $u =
O(D^3)$. For each $i=1,\ldots,u$, let $C_i$ denote the set of
curves of $C$ that intersect $\tau_i$, and let $P_i$ denote the set
of points that are contained in $\tau_i$. We set $m_i=|P_i|$ and
$n_i=|C_i|$, for $i=1,\ldots,u$, put $m' = \sum_i m_i \le m$, and
notice that, by construction, $n_i = O(n/D^2)$, for each $i$. 
An obvious property (which is a consequence of the generalized 
version of B\'ezout's theorem~\cite{Fu84}) is that each curve 
of $C$ intersects at most $ED+1 = O(D)$ cells of $\RR^3\setminus Z(f)$.

When $a'n^{1/k}\le m\le an^{3/2}$ (where the left inequality holds
by the induction assumption), we use, within each cell $\tau_i$ 
of the partition, for $i=1,\ldots, u$, the naive bound
$$
I(P_i,C_i) = O(m_in_i^{1-1/k} + n_i) = O\left(m_i(n/D^2)^{1-1/k}
+ n/D^2 \right) ,
$$
and, summing over the $O(D^3)$ cells, we get a total of
$$
O\left( \frac{ mn^{1-1/k} }{D^{2(1-1/k)}} + nD \right).
$$
With the above choice of $D$, we deduce that the total number of
incidences within the cells is
$$
O\left(m^{\frac{k}{3k-2}}n^{\frac{3k-3}{3k-2}}\right).
$$
When $m>an^{3/2}$, within each cell $\tau_i$ of the partition we
have $n_i = O(n/D^{2})=O(1)$, so the number of incidences within
$\tau_i$ is at most $O(m_in_i) = O(m_i)$, for a total of $O(m)$ incidences.
Putting these two alternative bounds together, we get a total of
\begin{equation}
\label{eq:th6}
O\left(m^{\frac{k}{3k-2}}n^{\frac{3k-3}{3k-2}}+m\right)
\end{equation}
incidences within the cells.

\smallskip

\noindent {\bf Incidences within the zero set $Z(f)$.}
It remains to bound incidences with points that lie on $Z(f)$.
Set $P^*:= P \cap Z(f)$ and $m^*:=|P^*|=m-m'$. Let $C^*$ denote
the set of curves that are fully contained in $Z(f)$, and set
$C': = C \setminus C^*$,
$n^* := |C^*|$, and $n' := |C'| =n-n^*$. Since every curve of
$C'$ intersects $Z(f)$ in at most $ED = O(D)$ points, we have (for
either choice of $D$)
\begin{equation} \label{eq:inc'}
I(P^*,C') = O(nD) =
O\left(m^{\frac{k}{3k-2}}n^{\frac{3k-3}{3k-2}}+m\right) .
\end{equation}
Finally, we consider the number of incidences between points of
$P^*$ and curves of $C^*$. Decompose $f$ into (complex)
irreducible factors $f_1,\ldots,f_t$, for $t\le \Deg$ (this order
of the factors is arbitrary but fixed), and assign
each point $p\in P^*$ (resp., curve $\gamma\in C^*$) to the first
irreducible factor $f_i$, such that $Z(f_i)$ contains $p$ (resp.,
fully contains $\gamma$; such a component always exists). The number 
of ``cross-incidences'', between points and curves assigned to
different factors, is easily seen, arguing as above, to be $O(nD)$,
which satisfies our bound. In what follows, we recycle the symbols $m_i$
(resp., $n_i$), to denote the number of points (resp., curves) assigned to
$f_i$, and put $D_i=\deg(f_i)$, for $i=1,\ldots, t$. We clearly have
$\sum_i m_i = |P^*| = m^*$, $\sum_i n_i = |C^*| = n^*$, and
$\sum_i D_i\le\deg(f)=\Deg$.

For each $i=1,\ldots, t$, there are two cases to consider.

\smallskip

\noindent{\bf Case 1: $Z(f_i)$ is infinitely ruled by curves of $\C_0$.}
By assumption, there are at most $q$ curves of $C$ on $Z(f_i)$, implying that $n_i \le q$. 
Put $m_{\rm inf} := \sum_i m_i$, summed over the infinitely ruled components $Z(f_i)$.
We project the points of $P_i$ and the curves of $C_i$ onto some generic
plane $\pi_0$. A suitable choice of $\pi_0$ guarantees that 
(i) no pair of intersection points or points of $P_i$ project to the same point, 
(ii) if $p$ is not incident to a curve $\gamma\in C_i$ then the projections of $p$ and of $\gamma$ remain non-incident, 
(iii) no pair of curves in $C_i$ have overlapping projections, and 
(iv) no curve of $C_i$ contains any segment orthogonal to $\pi_0$. 

Moreover, the number of degrees of freedom does not change in the projection,
if one uses a restricted version of this notion, defined shortly (see Sharir et al.~\cite{SSS}). 
Before discussing this, we note that this statement is false for the standard definition
of degrees of freedom. Consider for example the family $\C$ of circles in $\reals^3$,
whose projections onto some plane forms the family $\C^*$ of all ellipses in the plane.
In this case $\C$ has three degrees of freedom, but $\C^*$ has five degrees of freedom.
However, the property that the number of degrees of freedom does not increase in the 
projection holds if we define it only with respect to the projection $P^*$ of the given 
point set $P$ and the projection $C^*$ of the given set $C$ of curves, 
and assume that the projection plane is sufficiently generic. 
Specifically, we say that the projected set $C^*$ of $C$ has $k$ 
degrees of freedom \emph{with respect to $P^*$} if for any $k$-tuple of distinct 
points of $P^*$ there are at most $k$ curves of $C^*$ that pass through all of them.

We claim that if the original family $\C_0$ in $\reals^3$ has $k$ degrees of freedom 
then the projection $C^*$ of any finite set $C\subset\C_0$ does indeed have $k$ 
degrees of freedom with respect to the projection $P^*$ of any finite point set
$P\subset\reals^3$, assuming that the projection plane $\pi_0$ is sufficiently 
generic, as above. Indeed, let $p_1^*,\ldots,p_k^*$ be a $k$-tuple of distinct 
points of $P^*$, and let $\gamma^*\in C^*$ be the
projection of some curve $\gamma$ of $C$ that passes through all these points (if 
there is no such curve, there is nothing to argue about). For each $j=1,\ldots,k$, 
let $p_j$ be the unique point of $P$ that projects to $p_j^*$. By property (ii) of 
the genericity of $\pi_0$, $p_j\in\gamma$. Since there are at most $k$ curves of 
$C$ that pass through $p_1,\ldots,p_k$, the claim follows. See also \cite{SSS} 
for a related argument.

We remark that, as is easily checked, the above property also holds for the 
infinite family of all the curevs of $\C_0$ that are contained in some fixed 
algebraic surface, which is the only context in which we use this property 
(where the surface is an infinitely ruled component of $Z(f)$).

The number of incidences for the points and curves assigned to the same
$Z(f_i)$ is therefore equal to the number of incidences between the projected 
points and curves, which, by Theorem~\ref{th:PS}, is\footnote{%
  An inspection of Pach and Sharir~\cite{PS} shows that the analysis there
  too only requires that the given family of curves have $k$ degrees of freedom
  with respect to the given point set.}
$$
O\left(m_i^{\frac{k}{2k-1}}n_i^{\frac{2k-2}{2k-1}}+m_i+n_i\right) =
O\left(m_i^{\frac{k}{2k-1}}n_i^{\frac{k-1}{2k-1}}q^{\frac{k-1}{2k-1}}+m_i+n_i\right).
$$
Summing over $i=1,\ldots, t$, and using H\"older's inequality, we get the bound
$$
O\left(m^{\frac{k}{2k-1}}n^{\frac{k-1}{2k-1}}q^{\frac{k-1}{2k-1}}+m_{\rm inf}+n\right) ,
$$
which, by making $A$ sufficiently large, is at most
\begin{equation} \label{infrule}
\tfrac{A}{4}\left(m^{\frac{k}{2k-1}}n^{\frac{k-1}{2k-1}}q^{\frac{k-1}{2k-1}}+m_{\rm inf}+n\right) .
\end{equation}

\smallskip

\noindent{\bf Remark.}
This is the only step in the proof where
being of reduced dimension $s$, for $s$ sufficiently small, might
yield an improved bound (over the one in (\ref{infrule})); see
below, in the follow-up proof of Theorem~\ref{th:imprInc3}, where
this observation is exploited, for details.

\smallskip

\noindent{\bf Case 2: $Z(f_i)$ is not infinitely ruled by curves of
$\C_0$.} In this case, Theorem~\ref{th:ruledgz}(b) implies that
there exist suitable constants $c_1$, $t$ that depend on $E$ and on
the complexity of $\C_0$, such that there are at most $c_1 D_i^2$
exceptional curves, namely, curves that contain at least $c_1 D_i$
$t$-rich points. 
Therefore, by choosing $c$ (in the definition of $D$) sufficiently
small, we can ensure that, in both cases (of medium-range $m$ and large $m$), 
$\sum_i D_i^2 \le (\sum_i D_i)^2 = D^2\le n/(64c_1)$, say. 
This allows us to apply induction on the number of curves, to handle the
exceptional curves. 

Concretely, denote by $m_{\rm rich}$ (resp., $m_{\rm poor}$)
the number of $t$-rich (resp., $t$-poor) points assigned to components 
of $Z(f)$ that are not ruled by curves of $\C_0$. We thus have
$m_{\rm rich} + m_{\rm poor} + m_{\rm inf} = m^*$. Incidences
with the $t$-poor points are easy to bound, because each such point
is incident to at most $t$ curves of $C^*$ (that are assigned to the
same component as the point), for a total of $tm_{\rm poor} = O(m_{\rm poor})$
incidences. We thus continue the analysis with the $t$-rich points only.
For each $i$ (for which $Z(f_i)$ is not infinitely ruled by curves of $\C_0$),
we have an inductive instance of the problem involving at most $m_i$ 
$t$-rich points and at most $c_1 D_i^2 \le n/64$ curves of $C$. 
By the induction hypothesis, the corresponding incidence bound is at most
$$
A \left(m_i^{\frac{k}{3k-2}}(c_1 D_i^2)^{\frac{3k-3}{3k-2}} +
m_i^{\frac{k}{2k-1}}(c_1 D_i^2)^{\frac{k-1}{2k-1}}q^{\frac{k-1}{2k-1}}
+ m_i + c_1 D_i^2 \right) .
$$
We now sum over $i$. For the first terms, we bound each $m_i$ by $m$, and, 
for the fourth terms too, use the fact that $\sum_i D_i^\alpha \le D^\alpha$ 
for any $\alpha \ge 1$. For the third terms, we bound $\sum_i m_i$ by $m_{\rm rich}$.
For the second terms, we use H\"older's inequality. Overall, we get the incidence bound
\begin{gather*}
A \left( m^{\frac{k}{3k-2}}(c_1D^2)^{\frac{3k-3}{3k-2}}
+ m^{\frac{k}{2k-1}}(c_1D^2)^{\frac{k-1}{2k-1}}q^{\frac{k-1}{2k-1}}
+ m_{\rm rich} + c_1 D^2 \right) \\
\le A \left( m^{\frac{k}{3k-2}}(n/64)^{\frac{3k-3}{3k-2}}
+ m^{\frac{k}{2k-1}}(n/64)^{\frac{k-1}{2k-1}}q^{\frac{k-1}{2k-1}}
+ m_{\rm rich} + n/64 \right) ,
\end{gather*}
which can be upper bounded (since $k\ge 2$) by
\begin{equation} \label{eq:induct}
\tfrac{A}{4}\left(m^{\frac{k}{3k-2}} n^{\frac{3k-3}{3k-2}} +
m^{\frac{k}{2k-1}}n^{\frac{k-1}{2k-1}}q^{\frac{k-1}{2k-1}}+n \right) + Am_{\rm rich} .
\end{equation}

Except for these incidences, for each $f_i$, each non-exceptional
curve in $C^*$ that is assigned to $Z(f_i)$ is incident to at most
$c_1 D_i$ $t$-rich points; the total number of incidences of this kind 
involving the $n_i$ curves assigned to $Z(f_i)$ and their incident 
$t$-rich points is at most $n_i\cdot c_1 D_i = O(n_iD_i)$. 
Other incidences involving the non-exceptional curves
in $C$ that are assigned to $Z(f_i)$ only involve $t$-poor points 
that are assigned to $Z(f_i)$; as argued,
the overall number of such point-curve incidences is $O(m_{\rm poor})$.
Therefore, when $Z(f_i)$ is not infinitely ruled by
curves of $C^*$, the number of incidences with the $t$-rich points
assigned to $Z(f_i)$ is $O(n_iD_i)$, plus terms that are accounted 
for by the induction. Summing over these components $Z(f_i)$, we 
get the bound $O(nD)$ plus the inductive bounds in (\ref{eq:induct}), and, 
choosing $A$ to be sufficiently large, these bounds will collectively be at most
\begin{equation}
\tfrac{A}{2} \left(m^{\frac{k}{3k-2}}n^{\frac{3k-3}{3k-2}} +
m^{\frac{k}{2k-1}}n^{\frac{k-1}{2k-1}}q^{\frac{k-1}{2k-1}}+n\right) + A(m_{\rm rich} + m_{\rm poor}) . 
\end{equation}
Adding up all the bounds obtained so far, and choosing $A$ to be a 
sufficiently large constant, the number of incidences satisfies the 
inequality in~(\ref{eqind}), thus establishing the induction step,
and thereby completing the proof.
\proofend

\paragraph{Proof of Theorem~\ref{th:imprInc3}.}
The proof proceeds by the same double induction on $n$ and $m$, and
establishes the bound, for any prespecified $\eps>0$,
\begin{equation} \label{eqind-impr}
I(P,C) \le A m^{\frac{k}{3k-2}}n^{\frac{3k-3}{3k-2}} + A_\eps \left(
m^{\frac{2s}{5s-4}}n^{\frac{3s-4}{5s-4}}q^{\frac{2s-2}{5s-4}+\eps}+m^{2/3}n^{1/3}q^{1/3}+m+n\right) ,
\end{equation}
for a suitable constant $A$ that depends on $k$, $\mu$, $s$, $E$,
and the complexity of $\C_0$, and another constant $A_\eps$ that also depends on $\eps$. The flow of the proof is very similar
to that of the preceding proof. The main difference is in the case
where some component $Z(f_i)$ of $Z(f)$ is infinitely ruled by
curves from $\C_0$. Again, in this case it contains at most $q$ curves of $C^*$.

We take the points of $P^*$ and the curves of $C^*$ that are
assigned to $Z(f_i)$, and project them onto some generic plane
$\pi_0$ (the same, sufficiently generic plane can be used for all such components), 
as in the proof of Theorem~\ref{th:mainInc3}, and get the same properties (i)--(iv)
of the projected points and curves. Let $P_i$ and $C_i$ denote, respectively, the set 
of projected points and the set of projected curves; the latter is a set of $n_i$ plane
irreducible algebraic curves of constant maximum degree\footnote{%
  A projection preserves irreducibility and does not increase the degree;
  see, e.g., Harris~\cite{Harris} for a reference to these facts.}
$DE$. Moreover, as in the preceding proof, the contribution of $Z(f_i)$ to $I(P^*,C^*)$ 
is equal to the number $I(P_i,C_i)$ of incidences between $P_i$ and $C_i$. We can now apply
Theorem~\ref{incPtCu} to $P_i$ and $C_i$. To do so, we first note:
\begin {lemma} \label{sames-curves}
$C_i$ is contained in an $s$-dimensional family of curves.
\end {lemma}
\noindent{\bf Proof.}
Here it is more convenient to work over the complex field $\cplx$ (see the general remark
in the introduction).
Let $\Pi_0$ denote the projection of $\cplx^3$ 
onto $\pi_0$ (now regarded as a complex plane). Let $\C_0(f_i)$ denote the family of 
the curves of $\C_0$ that are contained in $Z(f_i)$, and let $\tilde{\C}_0(f_i)$
denote the family of their projections onto $\pi_0$ (under $\Pi_0$).
Define the mapping $\psi: \C_0(f_i) \to \tilde{\C}_0(f_i)$, by
$\psi(\gamma)=\Pi_0(\gamma)$, for $\gamma\in \C_0(f_i)$. 
By Green and Morrison~\cite{GM} 
$\C_0(f_i)$ and $\tilde{\C}_0(f_i)$ are algebraic varieties and $\psi$ is a
(surjective) morphism from $\C_0(f_i)$ to $\tilde{\C}_0(f_i)$. In
general, if $\psi: X \mapsto Y$ is a surjective morphism of
algebraic varieties, then the dimension of $X$ is at least as large
as the dimension of $Y$. Indeed, Definition 11.1~in Harris~\cite{Harris}
defines the dimension via such a morphism, provided that it is finite-to-one. 
Therefore, $\tilde{\C}_0(f_i)$ is of dimension at most
$\dim(\C_0(f_i))=s$, and the proof of the lemma, which can now be
transported back to the reals, is complete. $\Box$

Applying Theorem~\ref{incPtCu} to the projected points and curves, we conclude
that the number of incidences for the points and curves assigned to $Z(f_i)$ is at most
$$
B_\eps \left(m_i^{\frac{2s}{5s-4}} n_i^{\frac{5s-6}{5s-4}+\eps} +
m_i^{2/3}n_i^{2/3} + m_i + n_i \right) \le B_\eps \left(m_i^{\frac{2s}{5s-4}}
n_i^{\frac{3s-4}{5s-4}}q^{\frac{2s-2}{5s-4}+\eps} +
m_i^{2/3}n_i^{1/3}q^{1/3} + m_i + n_i \right) ,
$$
with a suitable constant of proportionality $B_\eps$ that depends on $s$ and on $\eps$. 
Summing this bound over all such components $Z(f_i)$, and applying H\"older's inequality 
(twice), we get the bound
$$
B'_\eps \left(m^{\frac{2s}{5s-4}}n^{\frac{3s-4}{5s-4}}q^{\frac{2s-2}{5s-4}+\eps}
+ m^{2/3}n^{1/3}q^{1/3} + m + n \right) ,
$$
for another constant $B'_\eps$ proportional to $B_\eps$. 
By making $A_\eps$ sufficiently large, this bound is at most
$$
\tfrac{A_\eps}{4}\left(m^{\frac{2s}{5s-4}}n^{\frac{3s-4}{5s-4}}q^{\frac{2s-2}{5s-4}+\eps}
+ m^{2/3}n^{1/3}q^{1/3} + m + n \right) .
$$
The rest of the proof proceeds as the previous proof, more or less
verbatim, except that we need a more careful (albeit
straightforward) separate handling of the leading term, multiplied
by $A$, and the other terms, multiplied by $A_\eps$. The induction
step then establishes the bound in (\ref{eqind-impr}) in much the
same way as above. \proofend

\smallskip

\noindent{\bf Remarks. (1)}
As already mentioned in the introduction, the ``lower-order'' terms
$$
O\left(m^{\frac{2s}{5s-4}}n^{\frac{3s-4}{5s-4}}q^{\frac{2s-2}{5s-4}+\eps}
+ m^{2/3}n^{1/3}q^{1/3} + m + n\right)
$$
in the bound are (almost) ``best possible'' in the following sense. If the
bound in Theorem~\ref{incPtCu} were optimal, or nearly optimal, in
the worst case, for points and curves of $\C_0$ that lie in a
constant-degree surface $V$ that is infinitely ruled by such curves,
the same would also hold for the lower-order terms in the bound in
Theorem~\ref{th:imprInc3}.\footnote{%
  Theorem~\ref{incPtCu} is formulated, and proved in~\cite{SZ}, 
  only for \emph{plane} curves. Nevertheless, it also holds for 
  curves contained in a variety $V$ of constant degree, simply by 
  projecting the points and curves onto some generic plane, as done 
  in other steps of our analysis.} 
This is shown by a simple packing argument, in which we take $n/q$ 
generic copies of $V$, and place on each of them $mq/n$ points and 
$q$ curves, so as to obtain
$$
\Omega\left((mq/n)^{\frac{2s}{5s-4}}q^{\frac{5s-6}{5s-4}} +
(mq/n)^{2/3}q^{2/3} +mq/n + q\right)
$$
incidences on each copy, for a total of
$$
(n/q)\cdot\Omega\left((mq/n)^{\frac{2s}{5s-4}}q^{\frac{5s-6}{5s-4}}
+ (mq/n)^{2/3}q^{2/3} + mq/n + q\right) =
\Omega\left(m^{\frac{2s}{5s-4}}n^{\frac{3s-4}{5s-4}}q^{\frac{2s-2}{5s-4}}
+ m^{2/3}n^{1/3}q^{1/3} +m + n\right)
$$
incidences. (This construction works when $m>n/q$. Otherwise, the
bound is linear, and clearly best possible. Also, we assume that the
lower bound does not involve the factor $q^\eps$, to simplify the
reasoning.) In particular, this remark applies to the case of points
and circles, as discussed in Theorem~\ref{ptcirc}.

\smallskip

\noindent{\bf (2)} There is an
additional step in the proof in which the fact that $\C_0$ is of
some constant (not necessarily reduced) dimension $s'$ could lead to
an improved bound. This is the base case $m = O(n^{1/k})$, where we
use the K{\H o}v\'ari-S\'os-Tur\'an theorem to obtain a linear bound
on $I(P,C)$. Instead, we can use the result of Fox et
al.~\cite[Corollary 2.3]{FPSSZ}, and the fact that the incidence
graph does not contain $K_{k,\mu+1}$ as a subgraph, to show that,
when $m= O(n^{1/s'})$, the number of incidences is linear. The
problem is that here we need to use the dimension $s'$ of the entire
$\C_0$, rather than the reduced dimension $s$ (which, as we recall,
applies only to subsets of curves of $\C_0$ that lie on a variety that is 
infinitely ruled by curves of $\C_0$). Typically, as already noted, $s$ is
larger than $k$ (generally twice as large as $k$), making this
bootstrapping bound inferior to what we have. Still, in cases where
$s'$ happens to be smaller than $k$, this would lead to a further
improved incidence bounds, in which the leading term is also smaller.

\paragraph{Rich points.}
Theorems \ref{th:mainInc3} and \ref{th:imprInc3} can easily be restated as bounding the
number of $t$-\emph{rich points} for a set $C$ of curves with $k$ degrees of freedom
(and or reduced dimension $s$) in $\reals^3$, when $t$ is at least some sufficiently
large constant. (Here richness is defined with respect to the entire set $C$.)
The case $t=2$ is treated in Guth and Zahl~\cite{GZ}, and
the same bound that they obtain holds for larger values of $t$ (albeit
without an explicit dependence on $t$), smaller than the threshold in the following corollary.
\begin{corollary} \label{co:rich}
(a) Let $C$ be a set of $n$ irreducible algebraic curves, taken from some
constructible family $\C_0$, of constant complexity, of irreducible curves of degree at
most $E$ and with $k$ degrees of freedom (with some multiplicity
$\mu$) in $\RR^3$, and assume that no surface that is infinitely
ruled by curves of $\C_0$, or, alternatively, by curves of degree at
most $E$, contains more than $q$ curves of $C$ (e.g., make this
assumption for all surfaces of degree at most $100E^2$). Then there
exists some constant $t_0$, depending on $k$ (and $\mu$) and on
$\C_0$, or, more generally, on $E$, such that, for any $t\ge t_0$,
the number of $t$-rich points, namely points that are incident to 
at least $t$ curves of $C$, is
$$
O\left( \frac{n^{3/2}}{t^{\frac{3k-2}{2k-2}}} +
\frac{nq}{t^{\frac{2k-1}{k-1}}} + \frac{n}{t} \right) ,
$$
where the constant of proportionality depends on $k$, $E$, $\mu$,
and the complexity of $\C_0$.

\smallskip

\noindent
(b) If $\C_0$ is also of reduced dimension $s$, the bound on the number of $t$-rich points becomes
$$
O\left( \frac{n^{3/2}}{t^{\frac{3k-2}{2k-2}}} +
\frac{nq^{\frac{2s-2}{3s-4}+\eps}}{t^{\frac{5s-4}{3s-4}}} + \frac{n}{t} \right) ,
$$
where the constant of proportionality now also depends on $s$ and
$\eps$. (Actually, the first and third terms come with a constant that is
independent of $\eps$.)
\end{corollary}
\noindent{\bf Proof.}
Denoting by $m_t$ the number of $t$-rich points, the corollary is obtained by combining the
upper bound in Theorem~\ref{th:mainInc3} or Theorem~\ref{th:imprInc3} with the lower bound $tm_t$.
\proofend

The bound in (b) is an improvement, for $s=k$, when $q>t^{k+\eps'}$, for another
arbitrarily small parameter $\eps'>0$, which is linear in the prespecified $\eps$.
(To be more precise, this is an improvement at all only when the second term dominates the bound.)

It would be interesting to close the gap, by obtaining a $t$-dependent bound also for
values of $t$ between $3$ and $t_0$.  It does not seem that the technique in
Guth and Zahl~\cite{GZ} extends to this setup.

\section{Incidences between points and circles and similar triangles in $\reals^3$}
\label{se:sim}

We first briefly discuss the fairly straightforward proof of Theorem~\ref{ptcirc}. 
As already mentioned in the introduction, we have $k=s=3$, for the case of circles, 
so we can apply Theorem~\ref{th:imprInc3} in the context of circles, and obtain the bound
$$
I(P,C) = O\left( m^{3/7}n^{6/7} + m^{2/3}n^{1/3}q^{1/3} +
m^{6/11}n^{5/11}q^{4/11+\eps} + m + n \right) ,
$$
for any $\eps>0$, where $q$ is the maximum number of the given
circles that are coplanar or cospherical. In fact, the extension of the
planar bound (\ref{eq:impbo}) to higher dimensions, due to Aronov et
al.~\cite{AKS}, asserts that, for any set $C$ of circles in any
dimension, we have
\begin{equation} \label{ptci}
I(P,C) = O\left(m^{2/3}n^{2/3} + m^{6/11}n^{9/11}\log^{2/11}(m^3/n) + m + n \right) ,
\end{equation}
which is slightly better than the general bound of Sharir and Zahl~\cite{SZ}
(given in Theorem~\ref{incPtCu}). If we use this bound, instead of that in
Theorem~\ref{incPtCu}, in the proof of Theorem~\ref{th:imprInc3} (specialized for the case of circles),
we get the slight improvement (in which the two constants of proportionality are now absolute)
$$I(P,C) = O\left( m^{3/7}n^{6/7} + m^{2/3}n^{1/3}q^{1/3} +
m^{6/11}n^{5/11}q^{4/11}\log^{2/11}(m^3/q) + m + n \right),$$ which establishes Theorem~\ref{ptcirc}.
\proofend

\paragraph{The number of similar triangles.}
Theorem~\ref{ptcirc} has the following interesting application. Let $P$ be a set of $n$ points in $\reals^3$, and let $\Delta=abc$ be a fixed
given triangle. The goal is to bound the number, denoted as $S_\Delta(P)$, of triangles
spanned by $P$ and similar to $\Delta$. The best known upper bound for $S_\Delta(P)$,
obtained by Agarwal et al.~\cite{AAPS}, is $O(n^{13/6})$, and the proof that establishes this bound in \cite{AAPS} is fairly
involved. Using Theorem~\ref{ptcirc}, we obtain the following simple and fairly straightforward
improvement.
\begin{theorem} \label{thm:simi}
$S_\Delta(P) = O(n^{15/7})$.
\end{theorem}
\noindent{\bf Proof.}
Following a standard strategy, fix a pair $p,q$ of points in $P$, and consider the locus
$\gamma_{pq}$ of all points $r$ such that the triangle $pqr$ is similar to $\Delta$
(when $p,q,r$ correspond to $a,b,c$, respectively).
Clearly, $\gamma_{pq}$ is a circle whose axis (line passing through the center of
$\gamma_{pq}$ and perpendicular to its supporting plane) passes through $p$ and $q$.
Moreover, there exist at most two (ordered) pairs $p,q$ and $p',q'$ for which
$\gamma_{pq} = \gamma_{p'q'}$. Let $C$ denote the set of all these circles (counted
without multiplicity). Then $S_\Delta(P)$ is at most two thirds of the number $I(P,C)$
of incidences between the $n$ points of $P$ and the $N=O(n^2)$ circles of $C$.

By Theorem~\ref{ptcirc} we thus have
$$
S_\Delta(P) = O\left( n^{3/7}(n^2)^{6/7} + n^{2/3}(n^2)^{1/3}q^{1/3} +
n^{6/11}(n^2)^{5/11}q^{4/11}\log^{2/11} n + n^2 \right) ,
$$
where $q$ is the maximum number of circles in $C$ that are either coplanar or cospherical.
That is, we have
\begin{equation} \label{sdp}
S_\Delta(P) = O\left( n^{15/7} + n^{4/3}q^{1/3} + n^{16/11}q^{4/11}\log^{2/11} n + n^2 \right) .
\end{equation}
We claim that $q=O(n)$. This is easy for coplanarity, because, for any fixed
plane $\pi$, each point $p\in P$ can generate at most one circle $\gamma_{pq}$
in $C$ that is contained in $\pi$. Indeed, the axis of such a circle is perpendicular
to $\pi$ and passes through $p$. This fixes the center of $\gamma_{pq}$ (on $\pi$), and it is
easily checked that the radius is also fixed. A similar argument holds for cospherical
circles. Here too, for a fixed sphere $\sigma$, each point $p\in P$ that is not the
center $o$ of $\sigma$ can generate at most two circles $\gamma_{pq}$ in $C$ that are
contained in $\sigma$. This is because the axis of such a circle must pass through $o$,
which fixes the axis $\lambda$. The circle is then the intersection of $\sigma$ with 
the cone with apex $p$ and axis $\lambda$ whose generators form angle $\angle bac$ 
with $\lambda$. Since this cone intersects $\sigma$ in at most two circles, the 
claim follows. For $p=o$ there are at most $n-1$ additional such circles.

Hence, plugging $q=O(n)$ into (\ref{sdp}), we get $S_\Delta(P) = O(n^{15/7})$, as asserted.
\proofend

We remark that Zahl's recent bound \cite{Zahl} on point-circle incidences in $\reals^3$, 
as reviewed in the introduction, yields the weaker bound $O(n^{12/5})$.

\section{\bf Discussion}

In this paper we have made significant progress on incidence
problems involving points and fairly general families of algebraic curves in three dimensions. 
The study in this paper raises several interesting open problems.

\smallskip
\noindent{\bf (i)} As remarked above, a challenging open problem
is to characterize all the surfaces that are infinitely ruled by
algebraic curves of degree at most $E$ (or by certain classes
thereof), extending the known characterizations for lines and
circles. A weaker, albeit still hard problem is to reduce the upper
bound $100E^2$ on the degree of such a surface, perhaps all the way
down to $E+1$, or at least to $O(E)$.

\smallskip
\noindent{\bf (ii)} It would also be interesting to find additional
applications of the results of this paper, like the one with an
improved bound on the number of similar triangles in $\reals^3$,
given in Section~\ref{se:sim}. One direction to look at is the
analysis of other repeated patterns in a point set, such as
higher-dimensional congruent or similar simplices, which can
sometimes be reduced to point-sphere incidence problems; see
\cite{AAPS,AgS}.

As already mentioned, the results of this paper have recently 
been used in Sharir et al.~\cite{ShZl} for bounding incidences
between points and curves with almost two degrees of freedom.

\noindent{\bf (iii)} A potentially weak issue in our analysis,
manifested in the proof of our main theorems, is that in order
to bound the number of incidences between points and curves on some
variety $V$ of constant degree, we project the points and curves on
some generic plane and use a suitable planar bound, from
Theorem~\ref{th:PS} or Theorem~\ref{incPtCu}, to bound the number of
incidences between the projected points and curves. It would be very
interesting if one could obtain an improved bound, exploiting the
fact that the points and curves lie on a variety $V$ in $\reals^3$,
under suitable (natural) assumptions on $V$. Note that we only need
to apply this argument for surfaces $V$ that are infinitely ruled
by the family $\C_0$ of the given curves. Perhaps this restricted
setup could aid in improving the analysis.

\smallskip
\noindent{\bf (iv)} Finally, it would be challenging to extend the
results of this paper to higher dimensions, extending the result
by the authors in \cite{SS4d}, obtained for the case of lines, to more
general families of curves.


\appendix

\section{On surfaces ruled by curves} \label{ap:ruled}

In this appendix we review, and sketch the proofs, of several tools
from algebraic geometry that are required in our analysis, the main one of which is Theorem~\ref{th:ruledgz}.
These tools are presented in Guth and Zahl~\cite[Section 6]{GZ}, but
we reproduce them here, in a somewhat sketchy form, for the convenience
of the reader and in the interest of completeness. We work over the complex
field $\cplx$, but the results here also apply to our setting over the real
numbers (see \cite{GZ,SS4d} and a preceding remark for discussions of this issue).

A subset of $\cplx^N$ described by some polynomial equalities and one
non-equality, of the form
$$
\{p \in \cplx^N \mid f_1(p)=0, \ldots, f_r(p)=0, g(p)\ne 0\},
\quad \text{for\;} f_1,\ldots,f_r, g \in \cplx[x_1,\ldots, x_N] ,
$$
is called \emph{locally closed}. We recall that the (geometric) degree of an algebraic
variety $V\subset \cplx^N$ is defined as the number of intersection points of $V$ with
the intersection of $N-\dim(V)$ hyperplanes in general position
(see, e.g., Harris~\cite[Definition 18.1]{Harris}).
Locally closed sets have the following property.
\begin {theorem}[B\'ezout's inequality; B\"urgisser et al.~\protect{\cite[Theorem 8.28]{BCS}}] \label{th:bein}
Let $V$ be a nonempty locally closed set in $\cplx^N$, and let $H_1,\ldots, H_r$ be algebraic
hypersurfaces in $\cplx^N$. Then
$$
\deg(V\cap H_1\cap \cdots \cap H_r) \le \deg(V)\cdot \deg(H_1)\cdots \deg(H_r) .
$$
\end {theorem}

A constructible set $C$ is easily seen to be a union of locally closed sets.
Moreover, one can decompose $C$ uniquely as the union of irreducible locally closed sets
(namely, sets that cannot be written as the union of two nonempty and distinct locally closed sets).
By~B\"urgisser et al.~\cite[Definition 8.23]{BCS}), the degree of $C$ is the
sum of the degrees of its irreducible locally closed components. Theorem~\ref{th:bein}
implies that when a constructible set
$C$ has complexity $O(1)$, its degree is also $O(1)$. We
also have the following corollaries.
\begin {corollary} \label{co:bez}
Let
$$
X = \{p \in \cplx^N \mid f_1(p)=0, \ldots, f_r(p)=0, g(p)\ne 0\},
\quad \text{for\;} f_1,\ldots,f_r, g \in \cplx[x_1,\ldots, x_N] .
$$
If $X$ contains more than $\deg(f_1)\cdots \deg(f_r)$ points, then $X$ is infinite.
\end {corollary}
\noindent{\bf Proof.}
Assume that $X$ is finite. Put $V=\{p\in \cplx^N \mid g(p)\ne 0\}$.
By B\'ezout's inequality (Theorem~\ref{th:bein}), we have
$$
\deg(X)\le \deg(V) \cdot \deg(f_1)\cdots \deg(f_r) = \deg(f_1)\cdots \deg(f_r) ,
$$
where the equality $\deg(V)=1$ follows by the definition of the degree of 
locally closed sets (see, e.g., B\"urgisser et al.~\cite[Definition 8.23]{BCS}). 
When $X$ is finite, i.e., zero-dimensional, its degree is equal to the number 
of points in it, counted with multiplicities. This implies that the number 
of points in $X$ is at most $\deg(f_1)\cdots \deg(f_r)$, contradicting the
assumption of the theorem. Therefore, $X$ is infinite. 
\proofend

As an immediate consequence, we also have:
\begin {corollary} \label{co:const}
Let $C\subset \cplx^N$ be a constructible set and write it as the union
of locally closed sets $\bigcup_{i=1}^t X_i$, where
$$
X_i = \{p \in \cplx^N \mid f^i_1(p)=0, \ldots, f^i_{r_i}(p)=0, g^i(p)\ne 0\},
\quad \text{for\;} f^i_1,\ldots,f^i_{r_i}, g^i \in \cplx[x_1,\ldots, x_N] .
$$
If $C$ contains more than $\sum_{i=1}^t \deg(f^i_1)\cdots \deg(f^i_{r_i})$ points, then $C$ is infinite.
\end {corollary}
For a constructible set $C$, let $d(C)$ denote the minimum of
$\sum_{i=1}^t \deg(f^i_1)\cdots \deg(f^i_{r_i})$, as in Corollary~\ref{co:const},
over all possible decompositions of $C$ as the union of locally closed sets.
By B\'ezout's inequality (Theorem~\ref{th:bein}), it follows that $\deg(C)\le d(C)$.
Corollary~\ref{co:const} implies that if $C$ contains more than $d(C)$ points, then it is infinite.

Following Guth and Zahl~\cite[Section 4]{GZ}, we call an algebraic
curve $\gamma\subset \cplx^3$ a \emph{complete intersection} if
$\gamma=Z(P,Q)$ for some pair of polynomials $P,Q$. We let
$\cplx[x,y,z]_{\le E}$ denote the space of complex trivariate
polynomials of degree at most $E$, and choose an identification of
$\cplx[x,y,z]_{\le E}$ with\footnote{%
  Here $\binom {E+3} 3$ is the maximum number of monomials of the
  polynomials that we consider. For obvious reasons, the actual
  representation should be in the complex projective space
  $\cplx \P^{\binom {E+3} 3}$, but we use the many-to-one representation in
  $\cplx^{\binom {E+3} 3}$ for convenience.}
$\cplx^{\binom {E+3} 3}$.  We use the variable $\alpha$ to denote an
element of $(\cplx[x,y,z]_{\le E})^2$, and write
$$
\alpha = (P_\alpha, Q_\alpha)\in \left(\cplx[x,y,z]_{\le E}\right)^2
= \left(\cplx^{\binom {E+3} 3}\right)^2.
$$
Given an irreducible curve $\gamma$, we associate with it a choice
of $\alpha \in \left(\cplx^{\binom {E+3} 3}\right)^2$ such that
$\gamma$ is contained in $Z(P_\alpha,Q_\alpha)$, and the latter is a
curve (one can show that such an $\alpha$ always exists; see
Guth and Zahl~\cite[Lemma 4.2]{GZ} and also Basu and Sombra~\cite{BS14}).
Let $x\in \gamma$ be a non-singular point\footnote{%
  Given an irreducible curve in $\reals^3$, a point $x\in \gamma$ is non-singular
  if there are polynomials $f_1, f_2$ that vanish on $\gamma$ such that
  $\nabla f_1(x)$ and $\nabla f_2(x)$ are linearly independent.}
of $\gamma$; we say that $\alpha$ is
associated to $\gamma$ at $x$, if $\alpha$ is associated to
$\gamma$, and $\nabla P_\alpha(x)$ and $\nabla Q_\alpha(x)$ are
linearly independent. We refer the reader to \cite[Definition 4.1 and Lemma 4.2]{GZ}
for details. This is analogous to the works of Guth and Katz~\cite{GK2} and
of Sharir and Solomon~\cite{SS4d} for the special cases of parameterizing lines in
three and four dimensions, respectively.

In what follows, we fix a constructible set $\C_0 \subset \C_{3,E}$ of
irreducible curves of degree at most $E$ in 3-dimensional space (recall
that the entire family $\C_{3,E}$ is constructible). Following
\cite[Section 9]{GZ}, we call a point $p\in Z(f)$, for a given polynomial
$f\in\cplx[x,y,z]$, a \emph{$(t, \C_0, r)$-flecnode}, if there are
at least $t$ curves $\gamma_1,\ldots, \gamma_t \in \C_0$, such that, for
each $i=1,\ldots, t$, (i) $\gamma_i$ is incident to $p$, (ii) $p$ is a
non-singular point of $\gamma_i$, and (iii) $\gamma_i$ osculates to $Z(f)$
to order $r$ at $p$. This is a generalization of the notion of a
\emph{flecnodal} point, due to Salmon~\cite[Chapter XVII, Section III]{salmon}
(see also ~\cite{GK2, SS4d} for details).

With all this machinery, we can now present a (sketchy) proof of
Theorem~\ref{th:ruledgz}. The theorem is stated in Section~\ref{sec:intro},
and we recall it here. It is adapted from Guth and Zahl~\cite[Corollary 10.2]{GZ},
serves as a generalization of the Cayley--Salmon theorem on surfaces ruled by lines
(see, e.g., Guth and Katz~\cite{GK2}), and is closely related to Theorem~\ref{th:gz}
(also due to Guth and Zahl~\cite{GZ}). We only consider here the first part of the theorem,
as the second part is an easy consequence.

\smallskip

\noindent{\bf Theorem \XX.}
{\it For given integer parameters $c,E$, there are constants $c_1=c_1(c,E)$,
$r=r(c,E)$, and $t=t(c,E)$, such that the following holds. Let $f$ be a complex
irreducible polynomial of degree $D \gg E$, and let $\C_0 \subset \C_{3,E}$ be
a constructible set of complexity at most $c$. If there exist at least $c_1 D^2$
curves of $\C_0$, such that each of them is contained in $Z(f)$ and contains at
least $c_1 D$ points on $Z(f)$ that are $(t,\C_0,r)$-flecnodes, then $Z(f)$ is
infinitely ruled by curves from $\C_0$. In particular, if $Z(f)$ is not infinitely
ruled by curves from $\C_0$ then, except for at most $c_1 D^2$ \emph{exceptional}
curves, every curve in $\C_0$ that is fully contained in $Z(f)$ is incident to at
most $c_1 D$ $t$-rich points (points that are incident to at least $t$ curves in 
$\C_0$ that are also fully contained in $Z(f)$).}

\smallskip

\noindent {\bf Proof.}
For the time being, let $r$ be arbitrary.
By Guth and Zahl~\cite[Lemma 8.3 and Equation (8.1)]{GZ}, since $f$ is irreducible,
there exist $r$ polynomials\footnote{%
  To say that $h_j$ is a ploynomial in $\alpha$ (and $p$) means that it is a polynomial in the
  $2\binom{E+3}{3}$ coefficients of the monomials of the two polynomials in the pair $\alpha$
  (and in the coordinates $(x,y,z)$ of $p$).}
$h_j(\alpha,p)\in \cplx[\alpha,x,y,z]$, for $j=1,\ldots, r$, of
degree at most $b_j$ in $\alpha$ (where $b_j$ is a constant
depending on $j$ and on $E$), and of degree $O(D)$ in $p=(x,y,z)$,
with the following property: let $\gamma$ be an irreducible curve,
let $p$ be a non-singular point of $\gamma$, and let $\alpha$ be
associated to $\gamma$ at $p$, then $\gamma$ osculates to $Z(f)$ to
order $r$ at $p$ if and only if $h_j(\alpha,p)=0$ for $j=1,\ldots,
r$. (These polynomials are suitable representations of the first $r$
terms of the Taylor expansion of $f$ at $p$ along $\gamma$;
see~\cite[Section 6.2]{GZ} for this definition, and
also~\cite{GK2,SS4d} for the special cases of lines in $\reals^3$
and $\reals^4$, respectively.)

Regarding $p$ as fixed, the system $h_j(\alpha,p)=0$, for $j=1,\ldots,r$, in
conjunction with the constructible condition that $\alpha \in \C_0$, defines
a constructible set $\C_p$. By definition, we have
$d(\C_p) \le \left(\prod_{j=1}^r b_j\right) \cdot d(\C_0)$, which is a constant that
depends only on $r$ and $E$. By Corollary~\ref{co:const}, $\C_p$ is either infinite or
contains at most $d(\C_p)=O(1)$ points. By Guth and Zahl~\cite[Corollary 12.1]{GZ},
there exist a Zariski open set $\O$, and a sufficiently large constant $r_0$, that depend
on $\C_0$ and $E$ (recall that $E$ is also assumed to be a constant, 
and see~\cite[Theorem 8.1]{GZ} for the way $r_0$ is obtained),
such that if $p\in\O$ is a $(t,\C_0, r)$-flecnode, with $r\ge r_0$,
there are at least $t$ curves that are incident to $p$ and are fully
contained in $Z(f)$. Since, by assumption, there are at least $c_1
D^2$ curves, each containing at least $c_1 D$ $(t, \C_0, r)$-flecnodes, 
it follows from~\cite[Proposition 10.2]{GZ} that\footnote{
For this proposition we require that $r$, and therefore $E$ too, be
considered as constants compared with $D$.}
there exists a Zariski open subset ${\O}$ of $Z(f)$, all of whose
points are $(t,\C_0,r)$-flecnodes. As noted above,~\cite[Corollary 12.1]{GZ} 
then implies that every point of ${\O}$ is incident to at
least $t$ curves of degree at most $E$ that are fully contained in
$Z(f)$. As observed above, when $t\ge \left(\prod_{j=1}^r b_j\right) \cdot d(\C_0)$, 
a constant depending only on $\C_0$ and $E$, $Z(f)$ is infinitely 
ruled on this Zariski open set. By a simple argument
(a variant of which is given in~\cite[Lemma 6.1]{SS4dv}), we can
conclude that $Z(f)$ is infinitely ruled by curves from
$\C_0$, thus completing the proof. \proofend

\end{document}